\newcommand{\Z}{{\mathbb Z}}

\documentclass{amsart}

\numberwithin{equation}{section}

\newtheorem{theo}{{\sc Theorem}}

\newtheorem{lem}[theo]{{\sc Lemma}}
\newtheorem{prop}[theo]{{\sc Proposition}}

\newcommand{\R}{{\mathbb R}}

\usepackage{epsfig}
\usepackage{amscd}
\usepackage{amsmath, amssymb, amsthm}
\usepackage{graphics}
\usepackage{color}
\usepackage{dsfont}
\newtheorem*{main-theorem}{Main Theorem}
\newtheorem{proposition}{Proposition}[section]
\newtheorem{theorem}{Theorem}
\newtheorem*{old-thm}{Theorem}

\newtheorem{corollary}[proposition]{Corollary}

\theoremstyle{definition}

\newtheorem*{remark}{Remark}
\numberwithin{equation}{section}

\def\11{\mathds{1}}

\def\Ci{{\mathcal C}^\infty}

\def\O{{\mathcal O}}

\def\phi{\varphi}

\def\be{\begin{eqnarray*}}
\def\ee{\end{eqnarray*}}
\def\ben{\begin{eqnarray}}
\def\een{\end{eqnarray}}

\def\lll{\left\langle}
\def\rrr{\right\rangle}

\def\L2R{L_{\text{Rest}}^2}

\def\tchi{\tilde{\chi}}

\newcommand{\ccal}{\mathcal{C}}
\newcommand{\dcal}{\mathcal{D}}
\newcommand{\ecal}{\mathcal{E}}
\newcommand{\gcal}{\mathcal{G}}
\newcommand{\jcal}{\mathcal{J}}

\newcommand{\pcal}{\mathcal{P}}

\newcommand{\ncal}{\mathcal{N}}

\begin{document}
\title[Quantum ergodic restriction for Cauchy Data]{Quantum ergodic restriction for Cauchy Data: Interior QUE and restricted QUE}
\author{Hans Christianson}
\address{Department of Mathematics, UNC-Chapel Hill \\ CB\#3250
  Phillips Hall \\ Chapel Hill, NC 27599}
\email{hans@math.unc.edu}

\author{John  A. Toth}
\address{Department of Mathematics and Statistics, McGill University, Montreal, CANADA }
 \email{jtoth@math.mcgill.ca} 

\author{Steve Zelditch}
\address{Department of Mathematics, Northwestern  University,
Evanston, IL 60208-2370, USA} \email{zelditch@math.northwestern.edu}


\begin{abstract}
We prove a quantum ergodic restriction theorem for the  Cauchy
data of a sequence of quantum ergodic eigenfunctions on a hypersurface $H$
of  a Riemannian manifold $(M, g)$.  The technique
of proof is to use a Rellich type identity to relate
quantum ergodicity of Cauchy data on $H$  to quantum
ergodicity of eigenfunctions  on the global  manifold $M$.  This has the interesting
consequence that if the eigenfunctions are quantum unique ergodic on
the global manifold $M$, then the  Cauchy data is
automatically quantum unique ergodic on $H$ with respect to operators
whose symbols vanish to order one  on the glancing set of unit tangential directions to $H$.

\end{abstract}

\maketitle


\section{Introduction}
\label{introduction}

This article is concerned with the QER (quantum ergodic restriction) problem for hypersurfaces in
compact Riemannian manifolds $(M, g)$.  
 We consider the
eigenvalue problem on $M$
$$\left\{ \begin{array}{l}  -\Delta_g \phi_j = \lambda_j^2 \phi_j, \;\;\; \langle \phi_j, \phi_k \rangle = \delta_{jk} \\
\\
B \phi_j = 0 \;\; \mbox{on}\;\; \partial M \end{array} \right.,$$
where $\langle f, g \rangle = \int_M f \bar{g} dV$ ($dV$ is the
volume form of the metric) and where $B$ is the boundary operator,
e.g. $B \phi = \phi|_{\partial M}$ in the Dirichlet case or $B
\phi = \partial_{\nu} \phi|_{\partial M}$ in the Neumann case. We
also allow $\partial M = \emptyset$.   We  introduce  the Planck constant 
$h_j = \lambda_j^{-1}$;  for notational simplicity
we often drop the subscript $j$. 
We then  denote the eigenfunctions in the orthonormal basis by  $\phi_h$  and the eigenvalues by $h^{-2}$, so that the eigenvalue problem takes the semi-classical form,
\begin{equation} \label{heq} \left\{ \begin{array}{l}  (-h^2\Delta_g-1) \phi_h = 0, \\
\\
B \phi_h = 0 \;\; \mbox{on}\;\; \partial M \end{array} \right., \end{equation}
where $B = I$ or $B = hD_\nu$ in the Dirichlet or Neumann cases
respectively.  Also, $\Delta_g$ denotes the negative Laplacian, e.g. on flat $\R^n$ it denotes $\sum_{j  =1}^n \frac{\partial^2}{\partial x_j^2}.$

Let $H \subset M$ be a smooth hypersurface which does not meet
$\partial M$ if $\partial M \neq \emptyset$.  The main result of this article  (Theorem
\ref{T:Lef}) is that the  semiclassical Cauchy data
 \begin{equation} \label{CD} CD(\phi_h)  := \{(\phi_h |_{H}, \;
h D_{\nu} \phi_h |_{H}) \}
\end{equation} of eigenfunctions is {\it always}
 quantum ergodic along any hypersurface $H \subset M$ if the
eigenfunctions are quantum ergodic on the global manifold $M$.  The proof is a generalization of the boundary
case where $H = \partial M$, which was proved  in \cite{HaZe} and in \cite{Bu}. Our proof is modeled
on that of \cite{Bu}, developing ideas of \cite{GL} (see also
\cite{CTZ-2} for an abstract microlocal approach).   This automatic QER 
property of Cauchy data stands in contrast to the conditional nature of the QER property
for the Dirichlet data alone, which requires an ``asymmetry" condition on $H$ with respect
to geodesics  \cite{TZ1,TZ2,DZ}.  Note that in the boundary case $H = \partial M$, the Dirichlet resp. Neumann  boundary
condition kills one of the two components of the Cauchy data, so that the Cauchy data QER theorem
appears the same as the QER theorem for Neumann data (resp. Dirichlet data) alone. 

As emphasized below the statement of Theorem \ref{T:Lef}, the deduction of quantum ergodicity of the Cauchy
data from quantum ergodicity of the sequence of ambient eigenfunctions holds for the full original sequence. 
Hence, if the original sequence is a complete orthonormal basis of eigenfunctions, i.e. if $\Delta_g$ is QUE  (quantum uniquely ergodic), then the Cauchy data of the full orthonormal basis  is also QUE
on $H$ for   any embedded
orientable separating  hypersurface $H$. We refer to this as the QUER property.  However, it is not necessarily the case that the Cauchy data is QUER for 
the full algebra of pseudo-differential operators. 
  In Corollary \ref{COR} it is proved
that QUE on  $(M, g)$ implies QUER with respect to the subalgebra 
of semiclassical pseudodifferential operators on $H$  whose symbols vanish to order 1 
along $S^* H$.  The restriction on the symbols arises  because  the passage from QUE in the ambient manifold to QUER
on the hypersurface involves multiplying the symbols by a certain factor  which vanishes
to order one on  $S^*H$, i.e. the unit directions (co-)tangent to $H$. Therefore, QUE in the ambient manifold does not imply
QUER for all pseudodifferential operators on $H$, and indeed the test operators    damp out the 
possible modes which concentrate microlocally on $H$.  We nevertheless refer to it as a QUER property because it holds for the entire sequence of eigenfunctions; there is no need to remove
a subsequence of density zero for the subalgebra limits.

To state the results precisely, we introduce some notation. 
We work with  the semiclassical calculus of pseudo-differential operators as in the references \cite{Bu,DZ,HaZe,TZ1,Zw};  see also \S \ref{Appendix} for background. On both $H$ and $M$ we fix  (Weyl) quantizations $a \to a^w$ of semi-classical symbols to
semi-classical pseudo-differential operators.
 When it is necessary to indicate  which manifold is involved, we either write $Op_H(a)$ for
pseudo-differential operators on $H$ or we use capital letters $A^w(x, h D)$ to indicate operators on $M$ and
small letters $a^w(y, h D_y)$ to indicate operators on $H$.

A sequence of functions $u_{h_j}$ on a manifold $M$ indexed by a sequence of Planck constants  is said to be quantum ergodic with limit measure $d\mu$ if
$$\langle a^w(x,h_j D_x) u_j, u_j \rangle \to \omega(a_0) : = \int_{T^* M} a_0 d\mu, $$
for all zeroth order semi-classical pseudo-differential operators, 
where $a_0$ is the principal symbol of $a^w(x,h_jD_x)$. In the classical case of compact Riemannian manifolds $(M, g)$ with
ergodic geodesic flow of \cite{Schnirelman,Zel1, CdV,ZZ}, the $\{u_j\}$ were a subsequence of density one of  an orthonormal basis of eigenfunctions
of $\Delta_g$ and $d\mu$ is normalized Liouville measure  on $S^* M$.  The functional
$a \to \langle a^w(x,h_jD_x) u_j, u_j \rangle$ is often referred to as a microlocal lift (or a Wigner distribution),
and the  limit measure or state $\omega(a)$ is often called
a quantum limit or a semi-classical defect measure.  Thus, we  use the term `quantum ergodic sequence' in this article
to mean a sequence of functions whose microlocal lifts have a unique weak* limit.  For the results
of this article, it is not necessary to assume that the geodesic flow is ergodic; we show that  that Cauchy data of a 
quantum ergodic  sequence of eigenfunctions  in the ambient space $M$ are quantum ergodic on $H$, no matter
what mechanism leads to quantum ergodicity of the  original sequence.


 We introduce a hypersurface  $H$, which we assume to  be
 orientable,
embedded, and  separating in the sense that
$$M \backslash H = M_+ \cup M_-
$$
where $M_{\pm}$ are domains with boundary in $M$.
 This is not a  restrictive assumption
since we can arrange that any hypersurface is part of the boundary
of a domain.

Given a quantization $a \to Op_H(a)$  of semi-classical symbols $ a \in S^{0}_{sc}(H)$ of order zero (see Appendix \ref{Appendix}) to semi-classical pseudo-differential operators on $L^2(H)$,
we  define the microlocal lifts of the
Neumann data as the  linear functionals on $a \in S^{0}_{sc}(H)$ given by 
$$\mu_h^N(a): = \int_{B^* H} a  \, d\Phi_h^N : = \langle Op_H(a) h D_{\nu} \phi_h
|_{H}, h D_{\nu} \phi_h |_{H}\rangle_{L^2(H)}.  $$ 
We
also  define the {\it renormalized microlocal lifts} of the Dirichlet
data by
$$\mu_h^D(a): = \int_{B^* H} a \, d\Phi_h^{RD} : = \langle Op_H(a) (1 +
h^2 \Delta_H) \phi_{h} |_H,  \phi_{h}|_H \rangle_{L^2(H)}.
$$
Finally, we  define the microlocal lift $d \Phi_h^{CD}$ of the Cauchy data  to be the sum
\begin{equation} \label{WIGCD} d \Phi_h^{CD} := d \Phi_h^N + d \Phi_h^{RD}. \end{equation}
Here, $h^2 \Delta_H$ denotes the negative
tangential Laplacian for the induced metric on $H$, so that
the operator $(1+h^2 \Delta_H)$ is characteristic precisely on the
glancing set $S^*H$ of $H$.  Intuitively,  we  have
renormalized the Dirichlet data by damping out the whispering gallery
components.

The distributions $\mu_h^N, \mu_h^D$  are asymptotically positive, but are not
normalized to have mass one and may tend to infinity.  
They  depend on the choice of quantization, but their possible weak* limits as $h \to 0$
do not, and the results of the article are valid for any choice of  quantization. We refer to \S \ref{Appendix} or to \cite{Zw} for background on semi-classical microlocal analysis.

 Our
first result is that the Cauchy data of a sequence of quantum ergodic
eigenfunctions restricted to $H$ is
automatically QER for semiclassical pseudodifferential operators with
symbols vanishing on the glancing set $S^*H$, i.e. that
$$ d \Phi_{h}^{CD}  \to \omega, $$
where $$\omega(a)  =  \frac{4}{\mu(S^*
 M)} \int_{B^*H} a_0(x', \xi') (1 - | \xi' |^2)^{1/2} d \sigma $$ is the limit state of Theorem \ref{T:Lef}.
This was  proved in a different way
in \cite{TZ1} in the case of piecewise smooth Euclidean domains. The assumption 
$H \cap
\partial M = \emptyset$  is for  simplicity of exposition and because the case $H = \partial M$ is already
known.

\begin{theorem} \label{thm1}
\label{T:Lef} Suppose $H \subset M$ is a smooth, codimension $1$
embedded orientable separating hypersurface  and assume $H \cap
\partial M = \emptyset$.
Assume that $\{\phi_h\}$ is a quantum ergodic sequence of eigenfunctions \eqref{heq}.  Then the sequence 
$\{d \Phi_{h}^{CD} \}$ \eqref{WIGCD} of microlocal lifts of the Cauchy data of $\phi_h$ is quantum ergodic on $H$ in the sense that for any
 $a \in S^0_{sc}(H),$
$$\begin{array}{l}
\lll Op_H(a) h D_\nu \phi_h |_{H} ,  h D_\nu \phi_h |_H \rrr_{L^2(H)}  + \lll Op_H(a) (1 +
 h^2 \Delta_H) \phi_{h} |_H, \phi_{h} |_H
\rrr_{L^2(H)} \\ \\ \rightarrow_{h \to 0^+} \frac{4}{\mu(S^*
 M)} \int_{B^*H} a_0(x', \xi') (1 - | \xi' |^2)^{1/2} d \sigma,
\end{array}$$
where $a_0(x', \xi')$ is the principal symbol of $Op_H(a)$, $- h^2
\Delta_H$ is the induced tangential (semiclassical) Laplacian with principal
symbol $| \xi' |^2$, $\mu$ is the Liouville measure on $S^*
M$, and $d \sigma$ is the standard symplectic volume form on $B^* H$.
\end{theorem}

\begin{remark}
We emphasize that the limit along $H$ in Theorem \ref{thm1} holds for the full sequence $\{\phi_h\}$.  Thus, if the full
sequence of eigenfunctions is known to be quantum ergodic, i.e. if the sequence is QUE,  then the conclusion of the theorem applies to the full sequence of
eigenfunctions.

We also remark that although we do not state it formally,  Theorem \ref{thm1} (and indeed all the
results in this paper) apply equally well to quasimodes, that is, to
approximate eigenfunctions satisfying
\[
\| (-h^2 \Delta_g -1) \phi_h \|_{L^2} = o(h) \| \phi_h \|_{L^2}
\]
as $h \to 0+$.

\end{remark}

The proof simply relates the interior and restricted microlocal lifts and reduces
the QER property along $H$ to the QE property of the ambient manifold.
 If we
assume that QUE holds in the ambient manifold, we automatically get
QUER, which is our first Corollary:

\begin{corollary} \label{COR} Suppose that $\{\phi_h\}$ is QUE on $M$.  Then
the  distributions $\{d \Phi_h^{CD}\}$ have a unique
weak* limit $$\omega(a): = \frac{4}{\mu(S^*
  M)} \int_{B^*H} a_0(x', \xi') (1 - | \xi' |^2)^{1/2} d \sigma.$$
  \end{corollary}

We note that $d \Phi_h^{CD}$ involves the microlocal lift $d \Phi_h^{RD}$ rather
than the microlocal lift of the Dirichlet data.  However, in Theorem
\ref{thm2}, we see that the analogue of  Theorem \ref{T:Lef} holds for a {\it density one
  subsequence} if we use the further renormalized distributions
$d \Phi_h^{D} + d \Phi_h^{RN}$ where
the microlocal lift $d \Phi_h^D \in
\dcal'(B^* H)$ of the Dirichlet data of $\phi_h$ is defined  by
$$\int_{B^* H} a \, d\Phi^D_h : = \langle Op_H(a) \phi_h|_{H}, \phi_h|_{H}
\rangle_{L^2(H)}, $$ and 
$$\int_{B^* H} a \, d\Phi^{RN}_h : = \langle (1 + h^2\Delta_H +
i0)^{-1} Op_H(a) h D_{\nu}\phi_h|_H, h D_{\nu}\phi_h |_{H} \rangle_{L^2(H)}. $$


\begin{theorem} \label{thm2}
Suppose $H \subset M$ is a smooth, codimension $1$
embedded orientable separating hypersurface  and assume $H \cap
\partial M = \emptyset$.
Assume that $\{\phi_h\}$ is a quantum ergodic sequence.  Then,  there exists a sub-sequence
of density one as $h \to 0^+$ such that for all  $a \in S^{0}_{sc}(H)$,
$$\begin{array}{l}
\lll (1 + h^2 \Delta_H + i0)^{-1}  Op_H(a)  h D_\nu \phi_h |_{H} , h D_\nu
\phi_h |_H \rrr_{L^2(H)} + \lll Op_H(a)  \phi_{h} |_H, \phi_{h} |_H
\rrr_{L^2(H)} \\ \\ \rightarrow_{h \to 0^+} \frac{4}{\mu(S^*
 M)} \int_{B^*H} a_0(x', \xi') (1 - | \xi' |^2)^{-1/2} d \sigma,
\end{array}$$
where $a_0(x', \xi')$ is the principal symbol of $Op_H(a).$ \end{theorem}


 The additional
step in the proof is  a pointwise  local Weyl law as in \cite{TZ1} section 8.4  showing that only a sparse
set of eigenfunctions could scar on the glancing  set $S^*H$. This is precisely 
the  step which is  not allowed in the QUER problem.  Therefore,
 QUER for all $Op_H(a)$ might fail for this rescaled problem; to determine whether it holds
for all $Op_H(a)$  we would need
a new idea.
However, the following is a direct consequence of Theorem \ref{thm2}.
\begin{corollary} \label{COR2} Suppose that $\{\phi_h\}$ is QUE on $M$.  Then
the  distributions $\{d \Phi_h^{D} + d \Phi_h^{RN}\}$ have a unique
weak* limit $$\omega(a): = \frac{4}{\mu(S^*
  M)} \int_{B^*H} a_0(x', \xi') (1 - | \xi' |^2)^{-1/2} d \sigma$$
   with respect to the subclass of symbols which vanish on
$S^*H$. \end{corollary}

We prove Theorem \ref{T:Lef}  by means of a Rellich identity adapted from \cite{GL,Bu}. 
It is also possible to prove the theorem using the layer potential approach in Step 2 (Proof
of (7.4)) in section 7 of \cite{HaZe}. To adapt this proof, one would need to introduce a semi-classical
Green's function in place of the Euclidean Green's function, verify that it has the properties
of the latter 
in section 4 of \cite{HaZe}, and then go through the proof of Step 2. Despite the authors'
fondness for the layer potential approach, this proof is much longer than the infinitesimal
Rellich identity approach and we have decided to omit the details.

\subsection*{Acknowledgements}
The first version of this article was written at the same time as \cite{TZ1,TZ2} but 
its completion  was post-poned while the authors proved the QER phenomenon for Dirichlet
data alone. We were further stimulated to complete the article 
by a discussion with  Peter Sarnak  at the Spectral Geometry
conference at Dartmouth in July, 2010 in which we debated whether QUE
in the ambient domain implies QUER along $H$. We said `yes', Sarnak said
`no';  Corollaries \ref{COR} and \ref{COR2} explain  the sense in which  both answers
are right.

The research of H.C. was partially supported by NSF grant \# DMS-0900524;
J.T. was partially supported by NSERC grant \# OGP0170280 and a
William Dawson Fellowship; S.Z. was partially supported by NSF grant  \# DMS-0904252.

\section{Rellich approach: Proof of Theorem \ref{T:Lef}}

We have assumed $H$ is a separating hypersurface, so that $H$ is 
the boundary of a smooth open submanifold
of $M$, $H = \partial M_+ \subset M$. There is no loss of generality in this assumption, 
since we may always use a cutoff to
a subset of $H$.   We then use a Rellich
type identity to write the integral of a commutator over
$M_+$ as a sum of integrals over the boundary (of course the same
argument would apply on $M_- = M \setminus M_+$).  The argument
is partially motivated by Burq's proof of boundary quantum
ergodicity (ie. the case $H = \partial M$).




 Let $x =(x_1,...,x_{n-1},x_n)= (x',x_n) $ be Fermi normal coordinates in a small
tubular neighbourhood $H(\epsilon)$ of $H $ defined near a point
$x_0 \in H$. In these coordinates we can locally write
$$H(\epsilon) := \{ (x',x_n) \in U \times {\mathbb R}, \, | x_{n} | < \epsilon \}.$$
Here $U \subset {\mathbb R}^{n-1}$ is a coordinate chart
containing $x_0 \in H$ and $\epsilon >0$ is arbitrarily small but
for the moment, fixed. We let $\chi \in C^{\infty}_{0}({\mathbb
R})$ be a cutoff with $\chi(x) = 0$ for $|x| \geq 1$ and $\chi(x)
= 1 $ for $|x| \leq 1/2.$
  In terms of the normal coordinates,
   $$-h^2\Delta_g = \frac{1}{g(x)} hD_{x_n} g(x) hD_{x_n}   + R(x_n,x',hD_{x'})$$
where, $R$ is a second-order $h$-differential operator along $H$
with coefficients that depend on $x_{n}$, and $R(0,x', hD_{x'}) = -h^2 \Delta_H$
is the
induced tangential semiclassical Laplacian on $H$.

Let $A(x, h D_x) \in \Psi_{sc}^{0}(M)$ be an order zero semiclassical pseudodifferential operator on $M$ (see \ref{Appendix}). By  Green's formula and \eqref{heq} we get  the Rellich identity
\begin{align}
\label{rellich}
  \frac{i}{h} \int_{M_+}  & \left( [-h^2 \Delta_{g},  \, A(x, h D_{x}) ] \, \phi_{h}(x)  \right) \overline{\phi_{h}(x)} \,  dx \\
= & \int_H \left(h D_{\nu}  \, A(x',x_n,h D_x)
 \phi_{h}|_{H} \right)  \overline{\phi_{h}}|_{H}  \,  d\sigma_{H} \notag \\
& +   \int_{H}      \left(\, A(x',x_n,h D_x) \,
  \phi_{h}|_{H} \right) \overline{hD_{\nu}
    \phi_{h}}|_{H}   \, d\sigma_{H} .\notag
\end{align}
Here, $D_{x_j} = \frac{1}{i} \frac{\partial}{\partial x_j}$,   $D_{x'}=(D_{x_1},...,D_{x_{n-1}}),$ \,  $D_{\nu} = \frac{1}{i} \partial_{\nu}$ where $\partial_{\nu}$
is the interior unit normal to $M_+$.

Given $a \in S^{0,0}(T^*H \times (0,h_0]),$ we then choose
$$ A(x',x_n,h D_x) = \chi ( \frac{x_n}{\epsilon}) \,h
D_{x_{n}}a^w(x', hD').$$
Since $\chi(0)=1$ it follows that the second term on the right side  of
(\ref{rellich}) is just


\begin{equation} \label{rellich2}
 \lll a^w(x', hD')h D_{x_n} \phi_h |_H, h D_{x_n} \phi_h |_H \rrr. \end{equation}
The  first term on right hand side   of (\ref{rellich}) equals
\begin{align} \label{rellich3}
  \int_{H} &    h D_{n} (\chi(x_n/\epsilon)h D_n a^w(x',
  hD')\phi_{h} )\Big|_{x_n=0}   \overline{ \phi_{h}}\Big|_{x_n=0}  \,
  d\sigma_{H} \\
= & \int_H \Big( \chi(x_n/\epsilon) a^w(x', hD') (h D_n)^2 \phi_h +
\frac{h}{i \epsilon} \chi'(x_n/\epsilon)hD_n a^w(x', hD') \phi_h
\Big)\Big|_{x_n = 0} \overline{\phi_h}\Big|_{x_n = 0} d \sigma_H \notag \\
= & \int_H ( \chi(x_n/\epsilon) a^w(x', hD') (1 - R(x_n, x', hD'))
\phi_h )\Big|_{x_n = 0} \overline{\phi_h}\Big|_{x_n = 0} d \sigma_H +
\O_\epsilon (h) \notag
,
\end{align}
since $\chi'(0) = 0$ and $((hD_n)^2 + R + O(h)) \phi_h = \phi_h$ in these
coordinates.  



It follows from (\ref{rellich})-(\ref{rellich3}) that
\begin{align} \label{rellich4}
&  \lll Op_H(a)  h D_\nu \phi_h |_{H} , h D_\nu \phi_h |_H \rrr_{L^2(H)} + \lll Op_H(a)  (1
+ h^2 \Delta_H) \phi_h |_{H}, \phi_h |_{H} \rrr_{L^2(H)}  \\
& =   \left\langle  \left( \left\{
       \xi_{n}^{2} +
      R(x_n,x',\xi') , \, \chi(\frac{x_n}{\epsilon})  \xi_{n}a(x', \xi')  \right\}  \right)^w\phi_{h}, \,\,
  \phi_{h} \right\rangle_{L^2(M_+)}+ {\mathcal
  O}_{\epsilon}(h).
\end{align}

We now assume that $\phi_h$ is a sequence of quantum ergodic eigenfunctions, and
 take the $h \rightarrow 0^+$ limit on both sides of
(\ref{rellich4}).  We  apply interior quantum ergodicity to the term on the right side of
\eqref{rellich4}.  We compute

\begin{equation} \label{COMM} \begin{array}{lll}  \left\{
       \xi_{n}^{2} +
      R(x_n,x',\xi') , \, \chi(\frac{x_n}{\epsilon})  \xi_{n}a(x',
      \xi')  \right\} & = &  \frac{2}{\epsilon} \chi'(\frac{x_n}{\epsilon})
    \xi_{n}^{2}a(x', \xi')\\ && \\ && + \chi(\frac{x_n}{\epsilon})  R_2(x',x_n,\xi') , \end{array}  \end{equation}
where $R_2$ is a zero order symbol. Let $\chi_2 \in \Ci$ satisfy
$\chi_2 (t) = 0$ for $t \leq -1/2$, $\chi_2(t) = 1$ for $t \geq
0$, and $\chi_2'(t) >0$ for $-1/2 < t < 0$, and let $\rho$ be a
boundary defining function for $M_+$. Then
$\chi_2(\rho/\delta)$ is $1$ on $M_+$ and $0$ outside a
$\delta/2$ neighbourhood.  Now the assumptions that the sequence
$\phi_h$ is quantum ergodic implies that the matrix element of the second term on the right side of \eqref{COMM} is bounded by
\begin{align*}
 \Big| \langle & ( \,
 \chi(x_n / \epsilon)  R_2(x, \xi') \,)^w
 \phi_{h}, \phi_{h} \rangle_{L^{2}(M_+)} \Big|\\
& \leq \| \chi_2 (\rho / \delta)   
 \chi(x_n / \epsilon) 
 \phi_{h} \|_{L^2(M)} \| \tchi_2(\rho / \delta) \tchi( x_n / \epsilon) \phi_{h} \|_{L^{2}(M)} \\
&  = {\mathcal O}_\delta(\epsilon) + o_{\delta,\epsilon}(1),
\end{align*}
where $\tchi$ and $\tchi_2$ are smooth, compactly supported functions
which are one on the support of $\chi$ and $\chi_2$ respectively.
Here, the last line follows from interior quantum ergodicity of the
$\phi_h$ since the volume of the supports of $\chi(x_n/\epsilon)$
and $\tchi(x_n / \epsilon)$ is comparable to
$\epsilon$.

To handle the matrix element of the first term on the right side of \eqref{COMM},  we note
that $\chi'(x_n/\epsilon)|_{M_+} = \tchi'(x_n/\epsilon)$ for
a smooth function $\tchi \in \Ci(M)$ satisfying $\tchi = 1$
in a neighbourhood of $M \setminus M_+$ and zero inside
a neighbourhood of $H$. Then, again by interior quantum
ergodicity, we have
\begin{align} \label{rellich5}
 & 2 \left\langle    \left( \frac{1}{\epsilon}
     \chi'(\frac{x_n}{\epsilon}) \, \xi_n^2 a(x',\xi') \right)^w \phi_{h}, \,\,
   \phi_{h} \right\rangle_{L^2(M_+)} \\
& = 2 \left\langle  \left( \frac{1}{\epsilon}
     \tchi'(\frac{x_n}{\epsilon}) \, \xi_n^2 a(x',\xi') \right)^w \phi_{h}, \,\,
   \phi_{h} \right\rangle_{L^2(M)} \notag \\
 &  =  \frac{2}{\mu(S^* M)}\int_{S^*M} \frac{1}{\epsilon}
 \tchi'(\frac{x_n}{\epsilon}) ( 1 - R(x',x_n,\xi') )a(x', \xi') \, d \mu +
 O(\epsilon) + o_{\epsilon}(1) \notag \\
& = \frac{2}{\mu(S^* M)}\int_{S^*M_+} \frac{1}{\epsilon}
 \chi'(\frac{x_n}{\epsilon}) ( 1 - R(x',x_n,\xi') ) a(x', \xi') \, d \mu +
 O(\epsilon) + o_{\epsilon}(1), \notag
\end{align}
since $\tchi'$ and $\chi'$ are supported inside $M_+$.
Combining the above calculations yields
 \begin{align} \label{rellich6}
 &  \lll Op_H(a)  h D_\nu \phi_h |_{H} , h D_\nu \phi_h |_H \rrr_{L^2(H)} + \lll Op_H(a)  (1 + h^2
\Delta_H) \phi_h |_{H},  \phi_h |_{H} \rrr_{L^2(H)}  \\
& =  \frac{2}{\mu(S^* M)}   \int_{S^*M_+}
\frac{1}{\epsilon} \chi'(\frac{x_n}{\epsilon}) ( 1 -
R(x',x_n,\xi') ) a(x', \xi') \, d \mu + O_\delta(\epsilon) +
o_{\delta,\epsilon}(1). \notag
\end{align}

Finally, we take the $h \rightarrow 0^+$-limit in
(\ref{rellich6})  followed by the $\epsilon \rightarrow
0^+$-limit, and finally the $\delta \to 0^+$ limit.  The result is
that, since the left-hand side in \eqref{rellich6} is independent
of $\epsilon$ and $\delta$,
\begin{align}
\lim_{h \rightarrow 0^+} & \lll Op_H(a)  h D_\nu \phi_h |_{H} , h
D_\nu \phi_h |_H \rrr_{L^2(H)} + \lll Op_H(a)  (1 + h^2 \Delta_H) \phi_h |_{H},
\phi_h |_{H} \rrr_{L^2(H)}  \notag \\
& = \frac{2}{\mu(S^* M)}  \int_{S^*_H M} ( 1 -
R(x',x_n=0,\xi')
) \, d \tilde{\sigma} \notag \\
& = \frac{4}{\mu(S^* M)}  \int_{B^* H} ( 1 - | \xi'|^2)^{1/2}
a(x', \xi') ) \, d {\sigma}, \label{upshot}
\end{align}
where $d \tilde{\sigma}$ is the symplectic volume form on $S^*_H
M$, and $d \sigma$ is the symplectic volume form on $B^* H$.

\qed


\section{Proof of Theorem \ref{thm2} and Corollary \ref{COR2}}
The proof follows as in Theorem \ref{thm1} with a few modifications. For fixed $\epsilon_1 >0$ we choose the test operator 
\begin{equation} \label{test2}
A(x',x_n,h D_x) :=  ( I + h^2 \Delta_H(x',h D') + i \epsilon_1 )^{-1} \chi ( \frac{x_n}{\epsilon}) \,h
D_{x_{n}}a^w(x', hD') \end{equation} 
and since $WF_{h}'( \phi_h|_{H}) \subset B^*H$ (see \cite{TZ2} section 11) it suffices to assume that $a \in C^{\infty}_{0}(T^*H)$ with
$$ \text{supp} \, a \subset  B^*_{1+\epsilon_1^2}(H).$$
Let $\chi_{\epsilon_1}(x',\xi') \in C^{\infty}_{0}( B^*_{1+
  \epsilon_1^2} (H)\setminus B^*_{1-2\epsilon_1^2} (H); [0,1])$ be a
cutoff near the glancing set $S^*H$ with $ \chi_{\epsilon_1} (x',\xi')
= 1$ when $(x',\xi') \in  B^*_{1+\epsilon_1^2} (H)\setminus
B^*_{1-\epsilon_1^2} (H).$ Then, with $A(x,h D_x)$ in (\ref{test2}), the same Rellich commutator argument as in Theorem \ref{thm1} gives
\begin{align} \label{thm2main} 
& \lll (1 + h^2 \Delta_H + i\epsilon_1)^{-1}  a^w(x',hD')
(1-\chi_{\epsilon_1})^w h D_\nu \phi_h |_{H} , h D_\nu \phi_h |_H \rrr_{L^2(H)} \\
& \quad + \lll a^w(x',hD')  (1-\chi_{\epsilon_1})^w \left( \frac{1 - | \xi'
    |^2}{1-|\xi'|^2 + i\epsilon_1}  \right)^w \phi_h |_{H}, \phi_h |_{H}
\rrr_{L^2(H)}  \notag \\ 
& \to \frac{4}{\mu(S^*
 M)} \int_{B^*H} a_0(x', \xi') (1-\chi_{\epsilon_1}(x',\xi')) \left(
 \frac{ (1 - | \xi' |^2)^{1/2} }{1-|\xi'|^2 + i\epsilon_1}  \right) \,
d \sigma. \notag 
\end{align}
It remains to determine the contribution of the glancing set  $S^*H$. As in
\cite{Bu,DZ,HaZe,TZ1} we use a local Weyl law to do this. Because of
the additional normal derivative term the argument is slightly
different than in the cited articles and so we give some details.
For the rest of this proof, we need to recall that $h \in \{
\lambda_j^{-1} \}$, and we write $h_j$ for this sequence to emphasize
that it is a discrete sequence of values $h_j \to 0$.  
 Since $\| a^w(x',hD') \|_{L^2 \rightarrow L^2} = O(1),$ it follows that for $h \in (0,h_0(\epsilon_1)]$ with $h_0(\epsilon_1) >0$ sufficiently small,
\begin{align} \label{weyl1}
& \frac{1}{N(h)} \sum_{h_j \geq h} \Big| \langle a^w(x',hD')  \chi_{\epsilon_1}^w
\phi_{h_j} |_{H}, \phi_{h_j} |_{H} \rangle_{L^2(H)}  \Big| \\
& \quad \leq C \frac{1}{N(h)} \sum_{h_j
  \geq h} \left( \| \chi_{\epsilon_1}^w \phi_{h_j} |_{H} \|_{L^2(H)} \,
\| \chi_{2\epsilon_1}^w  \phi_{h_j} |_{H} \|_{L^2(H)} +{\mathcal O}(h_j^{\infty})  \right) \notag \\ 
&   \quad \leq \frac{C}{2}
\frac{1}{N(h)} \sum_{h_j \geq h}     \left(   \|  \chi_{\epsilon_1}^w
  \phi_{h_j} |_{H} \|_{L^2(H)}^2 +  \| \chi_{2\epsilon_1}^w
  \phi_{h_j} |_{H} \|_{L^2(H)}^2 + {\mathcal O}(h_j^{\infty})
\right) \notag \\
& \quad = {\mathcal O}(\epsilon_1^2). \notag 
\end{align}

By a Fourier Tauberian argument  \cite{TZ1} section 8.4,  it follows that for $h \in (0,h_0(\epsilon_1)]$
\begin{equation} \label{weyl1.1}
 \frac{1}{N(h)} \sum_{h_j \geq h}       |  \chi_{\epsilon_1,2\epsilon_1}^w \phi_{h_j}|_H(x') |^2 
 = {\mathcal O}(\epsilon_1^2) \end{equation}
uniformly for $x' \in H.$ The last estimate in (\ref{weyl1}) follows from (\ref{weyl1.1}) by  integration  over $H.$


To estimate the normal derivative terms, we first recall the standard
resolvent estimate 
\[
 \|  ( 1 + h^2 \Delta_H + i\epsilon_1)^{-1} u \|_{H^2_h(H)} \leq
 C \epsilon_1^{-1} \| u \|_{L^2(H)},
\]
where $H^2_h$ is the semiclassical Sobolev space of order $2$ (see
\cite{Zw} Lemma 13.6).  
Applying the obvious embedding $H^2_h(H) \subset L^2(H)$, we recover  
\begin{align*}
 \|  ( 1 + h^2 \Delta_H + i\epsilon_1)^{-1} u \|_{L^2(H)} & \leq C \|
 ( 1 + h^2 \Delta_H + i\epsilon_1)^{-1} u \|_{H^2_h(H)} \\
& \leq   C\epsilon_1^{-1}\| u \|_{L^2(H)}
\end{align*}
to get that
\begin{align} \label{weyl2} 
& \frac{1}{N(h)} \sum_{h_j \geq h} | \langle (1 + h^2 \Delta_H +
i\epsilon_1)^{-1} a^w(x',hD') \chi_{\epsilon_1}^w  h_j D_{x_n} \phi_{h_j}|_{H},  h_j D_{x_n} \phi_{h_j}|_{H}
\rangle_{L^2(H)} | \\
& \quad \leq C' \epsilon_1^{-1} \frac{1}{N(h)} \sum_{h_j \geq h} \left(  \|  \chi_{\epsilon_1}^w hD_{x_n}\phi_{h_j}|_H \|_{L^2(H)}   \,  \| \chi_{2\epsilon_1}^w
hD_{x_n}\phi_{h_j}|_{H} \|_{L^2(H)}  + {\mathcal O}(h_j^{\infty}) \right) \notag \\ 
& \quad \leq \frac{C' \epsilon_1^{-1}}{2} \frac{1}{N(h)} \sum_{h_j \geq
  h}     \left(   \|  \chi_{\epsilon_1}^w hD_{x_n}\phi_{h_j}|_{H} \|_{L^2(H)}^2 +
  \| \chi_{2\epsilon_1}^w hD_{x_n} \phi_{h_j} |_{H} \|_{L^2(H)}^2 + {\mathcal O}(h_j^{\infty})
\right) \notag \\
& \quad =
{\mathcal O}(\epsilon_1^{-1} \epsilon_1^2) \notag \\
& \quad = {\mathcal
  O}(\epsilon_1). \notag 
\end{align}
The last estimate follows again from the Fourier Tauberian argument in \cite{TZ1} section 8.4, which gives
\begin{equation} \label{weyl2.1}
 \frac{1}{N(h)} \sum_{h_j \geq h}       |  \chi_{\epsilon_1,2\epsilon_1}^w h_j D_{x_n} \phi_{h_j}|_{H}(x') |^2 
 = {\mathcal O}(\epsilon_1^2) \end{equation}
 uniformly for $x' \in H.$

Since $\epsilon_1 >0$ is arbitrary, Theorem \ref{thm2} follows  from (\ref{weyl1}) and (\ref{weyl2}) 
 by letting $\epsilon_1 \rightarrow 0^+$ in (\ref{thm2main}).

\qed

\subsection{Proof of Corollary \ref{COR2}}

We now observe that Corollary \ref{COR2} follows almost immediately
from the proof of Theorem \ref{thm2}.  To see this, we notice that by
restricting our attention to symbols which vanish on the glancing set,
we do not need to pass through the local Weyl law/Tauberian argument,
which is the step by which one extracts a density one subsequence.
Hence the result applies to the full sequence.
\qed

\section{\label{Appendix} Appendix}

\subsection{ Semiclassical symbols} Let $M$ be a compact manifold. 
By a semiclassical symbol $ a \in S^{m,k}(T^*M \times [0,h_0))$, we mean a smooth function possessing an
asymptotic expansion as $h \to 0$ of the form,
\begin{equation} \label{sc1}
  a(x,\xi,h) \sim_{h \rightarrow 0^+} \sum_{j=0}^{\infty} a_{k-j}(x,\xi) h^{m+j}, \end{equation}
with $a_{k-j} \in S^{k}_{1,0}(T^*M) $.
Here, $S^{k}_{1,0}$ is the standard H\"{o}rmander class consisting of smooth functions $a(x,\xi)$ satisfying the estimates
$|\partial_{x}^{\alpha} \partial_{\xi}^{\beta} a(x,\xi)| \leq C_{\alpha,\beta} \langle \xi \rangle^{k-|\beta|}$ for all multi-indices $\alpha, \beta \in N^{n}.$ 
We say that $A(h) \in Op_h^w(S^{m,k}(T^*M \times [0,h_0))$ provided its Schwartz kernel is locally of the form
\begin{equation} \label{sc2}
A(h)(x,y) = (2\pi h)^{-n} \int_{\R^n} e^{i \langle x-y,\xi \rangle/h} a(\frac{x+y}{2},\xi,h) \, d\xi \end{equation}
with $a \in S^{m,k}.$  We denote the operator $A(h)$ by  $a^w(x,hD_x)$ (or simply $a^w$). 
By a symbol of order zero we mean that $a \in S^{0,0}$, and we refer to $a_0(x, \xi)$ as the principal symbol. In the latter case, we simply write 
$$S^{0}_{sc}(M):= S^{0,0}, \,\,\, \Psi_{sc}^0(M):=Op_h^w(S^{0,0}).$$ 
Finally, when $H \subset M$ is a hypersurface and $a \in S^{0}_{sc}(H),$ we sometimes write $Op_H(a) = a^w$ to indicate dependence on the submanifold, $H.$
We refer to \cite{Zw} for background.




\end{document}